\documentclass[12pt]{amsart} 
\usepackage{amsmath,amssymb,amscd}
\usepackage{graphics}

\newtheorem{theorem}{Theorem}[section]
\newtheorem{proposition}[theorem]{Proposition}
\newtheorem{definition}[theorem]{Definition}
\newtheorem{lemma}[theorem]{Lemma}
\newtheorem{corollary}[theorem]{Corollary}
\newtheorem{remark}[theorem]{Remark}
\newtheorem{example}[theorem]{Example}

\newcommand{\La}{\Lambda}
\newcommand{\Hn}{\hat{0}}
\newcommand{\He}{\hat{1}}
\newcommand{\A}{\mathcal{A}}
\newcommand{\Ablock}{\A^\ast}
\newcommand{\Abblock}{\A^{\ast \ast}}

\newcommand{\Kk}{{\bf k}}
\newcommand{\I}{\mathcal{I}}
\newcommand{\h}{\mathcal{H}}
\newcommand{\BAH}{\mathcal{B}_{\A,\h}}
\newcommand{\Sn}{\mathfrak{S}_n}
\newcommand{\Pn}{\mathfrak{P}_n}

\newcommand{\type}{\textsc{sh}}
\newcommand{\fib}{\textsc{fi}}

\newcommand{\Domn}{\mathrm{Dom}(n)}
\newcommand{\Refn}{\mathrm{Ref}(n)}

\newcommand{\Fq}{\mathbb{F}_q}
\newcommand{\GL}{\mathrm{GL}}

\begin{document}

\title{A note on blockers in posets}

\author[Bj\"orner and Hultman]{Anders Bj\"orner \and Axel Hultman}
\address {Department of Mathematics, Royal Institute of Technology,
SE-100 44, Stockholm, Sweden}
\email{bjorner,axel@math.kth.se} 


\date{March 4, 2004}

\begin{abstract}
The {\em blocker} $A^{*}$ of an antichain $A$ in a finite poset $P$ is the
set of elements minimal with the property of having with each member of $A$ 
a common predecessor.
The following is done: 
\begin{enumerate}
\item The posets $P$ for which $A^{**}=A$ for all antichains are characterized.
\item The blocker $A^*$ of a symmetric antichain in the partition
lattice is characterized.
\item Connections with the question of finding minimal size
blocking sets for certain set families
are discussed.
\end{enumerate}
\end{abstract}

\maketitle

\section{Introduction}
The blocker $A^*$ of a set family $A$ is a well-known construction in
combinatorics and combinatorial optimization.
Among the early references are \cite{lehman} and 
\cite{edmonds-fulkerson}, and the concept is discussed in several
elementary textbooks. A crucial property in this setting is that 
if $A$ is an antichain (no set contains another), then
$A^{**}=A$.

The construction of blockers can be directly generalized to antichains 
in any finite bounded
poset. In this paper we work in this generality. The generalized blocker
construction has previously been considered by Matveev \cite{matveev}
and by Bj\"orner, Peeva and Sidman \cite{BPS}.


For general posets all that remains of blocker duality is the relation
$A^{***}= A^*$, valid for every antichain $A$. The first question
we deal with is: {\em What posets have the property that
$A^{**}= A$ for {\em all} antichains $A$?} 
Such ``strong blocker duality'' is characterized in Section 2.

In \cite{BPS} symmetric antichains and their blockers
in the partition lattice $\Pi_n$ play
an important role due to their relevance for the theory of subspace
arrangements. The second question we address is:
{\em How does one compute the blocker of a symmetric
antichain in $\Pi_n$?} The answer, presented in Section 3,
involves both the dominance and the refinement orderings of
number partitions.

In the final section we discuss an algebraic approach to 
finding minimal size blocking sets to set families that can be 
realized as families of flats in a geometric lattice realizable
over a field. 

\section{Posets with strong blocker duality}\label{blocker}

We begin by agreeing on some notation.
A poset is {\em bounded} if it contains unique bottom and top
elements, denoted by $\Hn$ and $\He$, respectively.  Let $P$ be a
bounded poset. We denote by $\La$ its set of {\em atoms},
i.e. elements that cover $\Hn$, and given $x \in P$ we let $\La(x) \subseteq
\La$ be the set of atoms below $x$. 

If $P$ is a lattice, then $x \vee y$ and $x \wedge y$ denotes the
{\em join} (supremum) and {\em meet} (infimum), respectively, of two
elements $x, y \in P$.

We say that a set $A \subseteq P$ is an {\em antichain} if
$\hat{0}\notin A \neq \emptyset$, and the elements of $A$
are pairwise incomparable with respect to the partial ordering in $P$.

\begin{definition}
Let $A$ be an antichain in a finite bounded poset $P$. The {\em
  blocker} of $A$ is the antichain 
$$
A^*= \mathrm{\, min \,}  \, \{\, x\in P
\; | \;
\La(x) \cap \La(a) \neq \emptyset \,
\mathrm{ for \, every }\,  a\in A\,\}\, ,
$$
where $\mathrm{\, min \,} E$ denotes the set of minimal
elements of a subset $E\subseteq P$.
\end{definition}

\begin{remark}\label{boundedness}{\rm 
The requirement in this paper that $P$ is bounded is for convenience
only. The bottom 
element $\Hn$ plays no role whatsoever, and the top element $\He$
has as only function to make sure that $A^* \neq \emptyset$ for
all antichains $A$. Everything can be
reformulated for general (non-bounded) posets having at least one
element $x$ (not necessarily unique) 
above all its atoms. We have chosen the formulation
for bounded posets since this is notationally simpler, and
since the examples we have in mind are bounded.
}\end{remark}


A partial order
on the antichains in $P$ is
defined as follows:
we say that $A\le B$ for two antichains if for each $b\in B$ there
exists an $a\in A$ such that
$a\le b$.

\begin{lemma}[cf.\ \cite{BPS} and \cite{matveev}]\label{lem: duality} 
Let $A$ and $B$ be antichains in a finite bounded poset $P$. 
\begin{enumerate}
\item[(1)] If $A\le B$, then $B^*\le A^*$.
\item[(2)] $A^{**}\le A$.
\item[(3)] $A^{***}=A^{*}$
\end{enumerate}
\end{lemma}

\begin{proof}
The first two parts are straightforward from the definitions.
By part (2)
we get that
$A^{***}\le A^*$. On the other hand, part (1) applied to
$A^{**}\le A$ yields
$A^{***}\ge A^*$.
\end{proof}

\begin{remark}{\rm 
The poset of antichains in $P$ is in fact a distributive
lattice with meet operation \ $A\wedge B=  \mathrm{\, min \,} (A\cup B),$
on which the mapping $A \mapsto A^*$ is a Galois connection
(see e.g.\ \cite{aigner} for the definitions).
These properties are not used in what follows.
}\end{remark}

As was mentioned in the Introduction, 
the blocker construction is well-known for the special case when
$P={\mathbf 2}^V$ is the Boolean lattice of all subsets of a finite set $V$.
In this case,  $A^{**}=A$ for all antichains $A$.
This can be seen, for instance, by applying the next lemma to the case $P =
{\mathbf 2}^V$. 

\begin{lemma}\label{le:symmetry}
Let $V$ be a finite set, and suppose $P$ is an induced subposet of the
Boolean lattice ${\mathbf 2}^V$ such that $\emptyset$, $V$ and all
singletons belong to $P$. Then, for two antichains $A, B
\subset P$, we have $B = A^*$ if and only if the
following property is satisfied: 

{\sf Property $C$:} For all $U\in P$, $V\setminus U$ contains no
  member of $A$ if and only if $U$ contains a member of $B$.
\begin{proof}
Note that Property $C$ is equivalent to the assertion ``{\em for all $U
  \in P$, $U \cap a \neq \emptyset$ for all $a \in A$ if and only if $U
  \supseteq b$ for some $b \in B$}''. Thus, Property $C$ is
  satisfied if and only if $B$ is the antichain of minimal elements in
  the set $\{x \in P \mid x \cap a \neq \emptyset \text{ for all } a
  \in A\}$. This antichain is precisely $A^*$.
\end{proof}
\end{lemma}

\begin{definition} Let $V$ be a finite set. A subposet of the Boolean
  lattice ${\mathbf 2}^V$ induced by a family $\mathcal{S}\subseteq
  {\mathbf 2}^V$ is called {\em well-complemented}
  if {\rm (i)} the empty set and all
singletons belong to $\mathcal{S}$, and {\rm (ii)} $\mathcal{S}$ is closed
under taking complements in $V$.
\end{definition}

By the symmetry of Property $C$, it is immediate that $A^{**}=A$ for
all antichains $A$ in a well-complemented poset. In fact,
well-complemented posets are characterized by this property, as we now
show.

\begin{theorem}\label{th:bool}
Let $P$ be a finite bounded poset. Then the following are equivalent:
\begin{enumerate}
\item[(1)] $A^{**}=A$ for all antichains $A$ in $P$.
\item[(2)] $P$ is isomorphic to a well-complemented subposet of a
  Boolean lattice.
\item[(3)] $P$ satisfies
\begin{enumerate}
\item[(i)] if $\La(x)\subseteq \La(y) $ then
$x\le y$, for all $x,y\in P$,
\item[(ii)] for all $x\in P$ there exists $y\in P$ such that 
$\La\setminus \La(x)=\La(y).$
\end{enumerate}
\end{enumerate}
\end{theorem}
 
\begin{proof}
The implication (2) $\Rightarrow$ (1) follows from Lemma \ref{le:symmetry}.
\smallskip

We show that (1) $\Rightarrow$ (3).
Assume that $A^{**} = A$ for all antichains $A \subset P$. Note that,
in particular, this implies that the map $A\mapsto A^*$ is injective
on antichains in $P$. 
Suppose $\La(x) \subseteq \La(y)$ for some $x,y \in P$. If $x$ and $y$ are
incomparable, then $\{x, y\}^* = 
\{x\}^* = \La(x)$, contradicting injectivity of $A \mapsto A^*$. If,
instead, $x > y$, we must have $\La(x) = \La(y)$. A similar
contradiction is then obtained from $\{x\}^* = \{y\}^* = \La(x)$. We
conclude that $x\leq y$, proving part (i). 

Pick $x \in P$. We must show that $\La \setminus \La(x) = \La(y)$
for some $y \in P$, so suppose that this is not the case. Note that
$(\La \setminus 
\La(x))^* = \min \{z \in P \mid \La(z) \supset \La \setminus
\La(x)\}$. Hence, we have $\La(z) \cap \La(x) \neq \emptyset$ for all
$z \in (\La \setminus \La(x))^*$, implying that $x \geq t$ for some $t
\in (\La \setminus \La(x))^{**}$. This, however, implies $(\La
\setminus \La(x))^{**} \neq (\La \setminus \La(x))$, a contradiction.
\smallskip

It remains to show that (3) $\Rightarrow$ (2). Let  
${\mathbf 2}^{\La}$ denote the Boolean lattice of all subsets of the
set $\La$ of atoms in $P$. Define a map $\psi: P \to {\mathbf
  2}^{\La}$ by $x \mapsto \La(x)$. Clearly, $\psi$ is
order-preserving, and property (i) implies both that
$\psi$ is injective and that
the inverse mapping $\psi(P) \to P$
given by $\La(x) \mapsto x$ is order-preserving. Thus, $P$ is
isomorphic to $\psi(P)$. By construction, $\psi(P)$ contains all
singletons and the empty set, and property (ii) shows that
it is closed under taking complements.
\end{proof}

The theorem has the somewhat unexpected consequence that
strong blocker duality forces $P$ to be isomorphic to its
order dual.

\begin{corollary}
Suppose that $A^{**}=A$ for all antichains $A$ in $P$.
Then $P$ admits 
a fixed-point-free, order-reversing
bijection of order $2$ onto itself.
\end{corollary}
\begin{proof}
This is a direct consequence of the implication (1) $\Rightarrow$ (2).
\end{proof}

The equivalence (1) $\Leftrightarrow$ (2) shows that the posets with
strong blocker duality and $n$ labeled atoms are precisely the ones
obtained from the full Boolean lattice $ {\mathbf 2}^{\{1,\dots,
  n\}}$ by deleting an arbitrary 
family of complementary pairs of subsets, avoiding cardinalities
$0, 1, n-1, n$. Thus, there are 
$$2^{2^{n-1}-n-1}$$ such posets, and they are pairwise distinct. 
Dividing by the possible
symmetries we obtain the following estimate for the number $N_n$ of
nonisomorphic $n$-atom posets with strong blocker duality:
$$N_n\ge \frac{2^{2^{n-1}-n-1}}{n!}\ge\frac{2^{2^n /n}}{2^{n^{2}}}
= 2^{(2^n /n) -n^2}. 
$$
Out of this doubly-exponential number of posets there is, however,
only one that is a lattice.

\begin{corollary}
Let $L$ be a finite lattice. Then the following are equivalent:
\begin{enumerate}
\item[(1)] $A^{**}=A$ for all antichains $A$ in $L$.
\item[(2)] $L$ is Boolean.
\end{enumerate}
\end{corollary}
\begin{proof} We already know that (2) $\Rightarrow$ (1).
To prove $(1) \Rightarrow (2)$, it suffices to show that, for finite
$V$, the only well-complemented subposet of ${\mathbf 2}^V$ which is a
lattice is ${\mathbf 2}^V$ itself. Let $P \subset {\mathbf 2}^V$ be
another well-complemented subposet, and suppose $S\subset V$ is maximal
with the property $S \not \in P$. All coatoms (elements covered by
$\He = V$) and $V$ belong to $P$. Hence, $S$ is covered by more than one
element in ${\mathbf 2}^V$. This means that $\La(S)$ has multiple
minimal upper bounds in $P$, so that $P$ cannot be a lattice.


\end{proof}

As a small example, Figure 1 shows one of the three $4$-atom posets 
with strong blocker duality that are not lattices.
 
 \vspace{.2in} 

\begin{center}
\resizebox{!}{1.7in}{\includegraphics{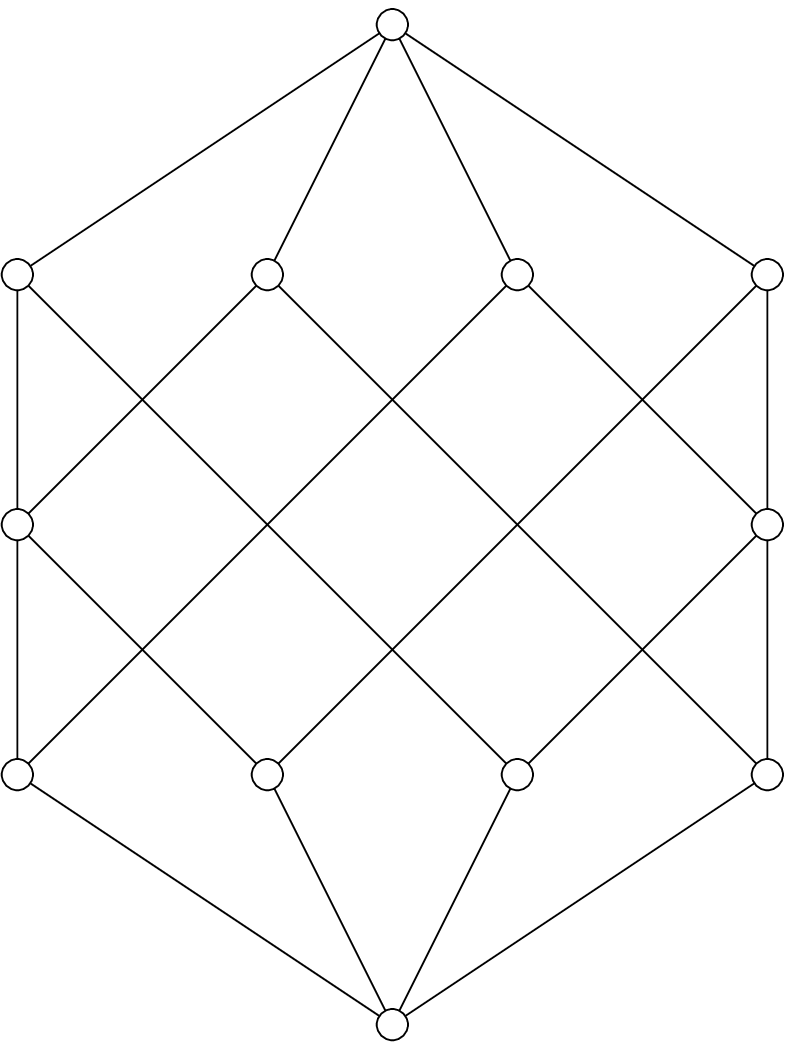}} \\
Figure 1.
\end{center}

\section{Symmetric blockers in partition lattices}
Recall that the {\em partition lattice} $\Pi_n$ consists of all set
partitions of $[n] = \{1,\dots,n\}$ ordered by refinement. In other
words, $\sigma \leq \tau \in \Pi_n$ if the equivalence relation
corresponding to $\tau$ contains the one corresponding to $\sigma$. 

We are interested in antichains in $\Pi_n$ that are invariant with
respect to the natural action of the symmetric group $\Sn$ on $\Pi_n$.
Since, clearly, the blocker of any $\Sn$-invariant antichain
is itself $\Sn$-invariant, the subject can be formulated
solely in terms of orbits, i.e.\ in terms of number partitions.

We need some notation. Let $\Pn$ be the set of 
partitions of the number $n$, $\Refn$ the refinement order on $\Pn$, and
$\Domn$ the dominance order. These partial orderings are
defined as follows. Let $\lambda=(\lambda_1, \lambda_2, \dots)$,
$\mu=(\mu_1, \mu_2, \dots)  \in \Pn$, with the parts 
$\lambda_i$ and $\mu_j$ decreasingly arranged and $\sum \lambda_i
=\sum \mu_j =n$. Then
\begin{enumerate}
\item $\lambda\le \mu$ in $\Domn$ if  $\sum_{i\le k} \lambda_i
\le\sum_{j\le k} \mu_j$ for all $k$,
\item $\lambda\le \mu$ in $\Refn$ if $\lambda$ can be obtained from $\mu$
by partitioning the parts $\mu_j$ .
\end{enumerate}
Note that the identity mapping $\Refn \rightarrow \Domn$ is order-preserving.

Let $\type: \Pi_n \to \Pn$ be the shape map $\{\tau_1, \dots, \tau_t\}
\mapsto \{|\tau_1|, \dots, |\tau_t|\}$ (multiset). In the other
direction, given $\lambda \in \Pn$ we let $\fib(\lambda)$ denote the 
fiber (inverse image) in $\Pi_n$, i.e.\ $\tau \in \fib(\lambda)$
iff $\type(\tau) = \lambda$. Similarly, for $S
\subseteq \Pn$ we define $\fib(S) = \type^{-1}(S)$.

The following theorem, characterizing blocker duality of
symmetric antichains in the partition lattice $\Pi_n$,
is based on the fact that every symmetric (i.e., $\Sn$-invariant)
 antichain in $\Pi_n$ is of the form $\fib(A)$ for some
 antichain $A$ in $\Refn$.

For a poset $P$ and $x \in P$, we write
$P_{\leq x} = \{y \in P \mid y \leq x\}$, and $P_{< x}$ is defined
similarly. The transpose of a number partition $\lambda$ is
denoted by $\lambda^\prime$. 

\begin{theorem} \label{th:blocker}
Let $A \subseteq \Refn$ be an antichain. We have
$\fib(A)^* = \fib(B)$ in $\Pi_n$, where 
\[
B = \min_{\Refn} \Pn \setminus \big(\cup_{\lambda \in A}\,
\Domn_{\leq \lambda^\prime}\big).
\]
In other words, to construct $B$ we take the refinement-minimal
number partitions among those that are not dominated by any $\lambda^\prime$,
$\lambda \in A$.
\begin{proof}
It suffices to show that, given $\lambda, \mu \in \Pn$, there exist
set partitions $\sigma \in \fib(\lambda)$ and $\tau \in \fib(\mu)$ with
$\sigma \wedge \tau = \Hn$ if and only if $\lambda^\prime$ dominates $\mu$.

We have a 1-1 correspondence between pairs of set partitions
$\sigma, \tau \in \Pi_n$ and bipartite graphs with $n$ labeled
edges and no isolated vertices as follows. Given $\sigma, \tau \in
\Pi_n$, the vertex set of the graph can be thought of as the set of
blocks in $\sigma$ and $\tau$. The graph is constructed by letting the
$i$-th edge connect the block containing $i$ in 
$\sigma$ and the block containing $i$ in $\tau$. The crucial
observation is that this graph
contains multiple edges if and only if $\tau \wedge \sigma \neq
\Hn$. By the Gale-Ryser Theorem, there is a bipartite graph with
degree sequences $\lambda, \mu \in \Pn$ without multiple edges if
and only if $\lambda^\prime$ dominates $\mu$. Hence the theorem.
\end{proof}
\end{theorem}

\begin{example}\label{hook}{\rm 
As a special case of the theorem we observe that 
the antichain $\fib((p,1^{n-p}))$ of all set partitions of hook 
type $(p,1^{n-p})$ is the blocker in $\Pi_n$ 
of the antichain $A_p$ consisting of all partitions with 
$p-1$ blocks. Conversely, $A_p$ is the blocker of
$\fib((p,1^{n-p}))$.}
\end{example}

\begin{corollary} \label{co:hook}
The blocker of every $\Sn$-invariant antichain in $\Pi_n$ contains a
hook shape antichain $\fib((p,1^{n-p}))$ for some $p$. In particular,
$\fib(\lambda)$ is itself  a blocker if and only if $\lambda$ is a
hook shape.
\begin{proof}
Let $A\subseteq \Refn$ be an antichain. In $\Refn$ as
well as in $\Domn$, the hook
shapes form a chain from the bottom element to the top element. Thus,
$\Pn \setminus (\cup_{\lambda \in A}\,
\Domn_{\leq \lambda^\prime})$ contains a unique $\mu =
(p,1^{n-p})$ which is both refinement-minimal and dominance-minimal among
the hook shapes. Now, $\Refn_{< \mu}$ has a unique dominance-maximal
element, namely $(p-1,1^{n-p+1})$. Thus, $\mu$ is
refinement-minimal in $\Pn \setminus (\cup_{\lambda \in A}\,
\Domn_{\leq \lambda^\prime})$.

For the last assertion, see Example \ref{hook}.
\end{proof} 
\end{corollary}

\begin{corollary}
The map $A\mapsto \fib(A)^*$ determines a 
bijection between antichains in $\Domn$
and $\Sn$-invariant blockers in $\Pi_n$. 

\begin{proof}
In this proof, let $X = \{A \subseteq \Refn \mid
A \text{ is an antichain}\}$ and $Y = \{A \subseteq
\Domn \mid A 
\text{ is an antichain}\}$. Define $\phi:X\to Y$ by letting $\phi(A)$
be the set of dominance-minimal elements in $A$. From the 
fact that the identity mapping
$\Refn \to \Domn$ is order-preserving follows that
$\phi$ is surjective.  

Note that, for $A, B\in X$, we have $\cup_{\lambda\in
  A}\Domn_{\leq \lambda^\prime} = \cup_{\lambda \in
  B}\Domn_{\leq \lambda^\prime}$ if and only if $\phi(A) =
  \phi(B)$. By Theorem \ref{th:blocker}, this implies $\phi(A) =
  \phi(B) \Leftrightarrow \fib(A)^* = \fib(B)^*$. Thus, $\phi(A)
  \mapsto \fib(\phi(A))^* = \fib(A)^*$ is a bijection from $Y$ to the set of
  $\Sn$-invariant blockers in $\Pi_n$.

\end{proof}
\end{corollary}

\section{Subspace arrangements and blocker ideals}

Here we review some necessary background for the following section.
This concerns subspace arrangements,
which provided the motivation for the blocker construction
in \cite{BPS}. 
For background and details concerning subspace arrangements,
see \cite{bjorner}.

Let $\Kk$ be a
field, and consider an arrangement $\A$ of subspaces of $\Kk^n$. The
{\em vanishing ideal} $\I_\A \subseteq \Kk[x_1,\dots,x_n]$ is the ideal
of polynomials that are identically zero on all subspaces in $\A$. 
It is an intriguing problem to determine generators
for $\I_\A$. 

Now, the arrangement $\A$ can always be embedded in a hyperplane
arrangement $\h$. In particular, $\A$ can be considered an antichain
in the intersection lattice $L_\h$  (a geometric lattice). In this
setting, we may define the {\em blocker ideal}
\[
\BAH = \langle \{\prod_{H \in \Lambda(B)}\ell_H \mid B \in \Ablock\} \rangle,
\]  
where $\ell_H$ is the defining linear form of the hyperplane $H$. 

It is easy to see that $\BAH \subseteq \I_\A$, and this inclusion is
in general strict. However, it turns out that in several of the cases
where generators for $\I_\A$ are known, we actually have $\BAH =
\I_\A$.

One particularly interesting and rich class of subspace arrangements
is the class of orbit arrangements, which we now define. The {\em
  braid arrangement} $\A_n$ is the arrangement of hyperplanes defined
by the equations $x_i = x_j$ for $1 \leq i < j \leq n$. Its
intersection lattice $L_{\A_n}$ is naturally isomorphic to the
partition lattice $\Pi_n$. The symmetric group $\Sn$ acts on the braid
arrangement by permuting the indices, and the subspace arrangements that
correspond to $\Sn$-invariant antichains in $L_{\A_n}$ we call {\em orbit
  arrangements}. As in the previous section, there is a 1-1
correspondence between orbit 
arrangements and antichains in $\Refn$. We let
  $\A_\lambda$ denote the arrangement corresponding 
to the partition $\lambda$.

Two interesting cases where it is known that $\BAH = \I_\A$ are when
$\A = \A_{(p,1^{n-p})}$ and $\A = \cup \A_\lambda$ (union over all
$\lambda$ with $p-1$ parts). These results are due  
to Li and Li \cite{li-li} and to Kleitman and Lov\'asz
\cite{lovasz}, respectively. In view of Example \ref{hook}, note that
(given $p$) either of the two arrangements is the blocker of the
  other. Actually, only blockers can be expected to have the
  property $\BAH = \I_\A$. This is so because of the following consequence of
  \cite[Theorem 3.3.4]{BPS}, 
if $\Kk$ is algebraically closed:
\begin{equation}
\BAH = \I_\A \quad\Longrightarrow\quad \Abblock = \A.
\end{equation}

\section{Minimal blocking sets}
Again, suppose $A$ is an antichain in a finite bounded poset $P$. We say
that a subset $S \subseteq \Lambda$ of the atoms is {\em
  $A$-intersecting} if $ S \cap \Lambda(a) \neq \emptyset$ for all $a
\in A$. Clearly, $\Lambda(b)$ is $A$-intersecting for every $b \in
A^*$. 

\begin{definition}
The antichain $A$ has the {\em Tur\'an property} if the smallest
cardinality of any $A$-intersecting atom set is $\min\{|\Lambda(b)|
\mid b \in A^*\}$. 
\end{definition}

To motivate this definition, again consider the antichain
$\fib((p,1^{n-p}))$ in
$\Pi_n$ consisting of all set partitions of shape $(p,1^{n-p})$ for
some fixed $p$. We may think of $\Pi_n$ as the lattice of all clique
graphs (i.e., graphs such that every connected component is a clique) 
on vertex set $[n]$, the atoms of $\Pi_n$
corresponding to the set of edges. Then,
the assertion that $\fib((p,1^{n-p}))$ has the Tur\'an property is
equivalent to the assertion that the smallest number of edges in any graph that
intersects every $p$-clique is attained in a clique graph on $p-1$
cliques. By passing to complements, one sees that
this is precisely the famous Tur\'an theorem of graph theory.

It seems reasonable to inquire which antichains have the Tur\'an
property. In particular, if antichains in $\Pi_n$ corresponding to
$\Sn$-orbits have the Tur\'an property, this gives rise to Tur\'an
type graph theorems.

In their paper, Li and Li \cite{li-li} point out that their theorem
implies the original Tur\'an theorem. Their argument can be
generalized to obtain the following.


\begin{theorem} \label{pr:turanideal}
Let $\A$ be a subspace arrangement embedded in a hyperplane
arrangement $\h$. If $\BAH = \I_\A$, then $\A$ has the Tur\'an
property (viewed as an antichain in $L_\h$).
\begin{proof}
Suppose $\A$ does not have the Tur\'an property. Then there exists a
set of hyperplanes $S \subseteq \h$ whose union contains all subspaces
in $\A$, and $|S| < |\Lambda(B)|$ for all $B \in \Ablock$. Thus, by
definition, we have $\deg(p) > |S|$ for all $p \in \BAH$. However, it
is easy to see that 
\[
\prod_{H \in S}\ell_H \in \I_\A,
\]
where, again, $\ell_H$ is the defining linear form of a hyperplane $H$. This
polynomial has degree $|S|$, and therefore $\BAH \neq \I_\A$. 
\end{proof}
\end{theorem}

\begin{example}{\rm 
We illustrate what this says with a small example,
where $\h$ is taken to be the braid arrangement $\A_6$
and hence $L_\h \cong \Pi_6$.

Let $A=\{222, 3111\}$ and $B=\{42, 51\}$ be two antichains
in $\mathrm{Ref}(6)$. One sees from Theorem \ref{th:blocker}
that $\fib(A)^* =\fib(B)$ and $\fib(B)^* =\fib(A)$ in $\Pi_6$. It
was checked in Example 3.4.3 of \cite{BPS} that the blocker ideal
equals the vanishing ideal for the corresponding orbit arrangements
in both cases. Thus, Theorem \ref{pr:turanideal} applies.

What the Tur\'an property then means in the case
$A=\{222, 3111\}$ is the following: {\em The maximal number of edges
of a graph on $6$ vertices not containing three independent edges or
a $3$-clique equals} \\
max$\{\#K_{4,2}, \#K_{5,1}\}$=max$\{8,5\}=8$.
Here $\#K_{n,m}$ denotes the number of edges in the complete
bipartite graph $K_{n,m}$. Note that if we excluded only a
$3$-clique, the answer would be 
max$\{\#K_{3,3}, \#K_{4,2}, \#K_{5,1}\}=9$,
which of course agrees with Tur\'an's theorem.

Similarly, what the Tur\'an property means in the case
$B=\{42, 51\}$ is: {\em The maximal number of edges
of a graph on $6$ vertices not containing either
a $4$-clique and an independent edge or
a $5$-clique, equals} \\
max$\{\#K_{2,2,2}, \#K_{3,1,1,1}\}$=max$\{12,12\}=12$.}
\end{example}

\begin{remark}{\rm 
The converse of Theorem \ref{pr:turanideal} does not hold in
general, not even for blockers. A construction of D.\ Kozlov (see
\cite[Example 4.2.2]{BPS}) yields an arrangement $\A$ of two subspaces
embedded in an arrangement $\h$ of four hyperplanes such
that $\A$ has the Tur\'an property, but $\BAH \neq \I_\A$. Moreover,
$\A$ is a blocker in $L_\h$.} 
\end{remark}

\begin{example}{\rm 
 A graph-theoretic theorem by Simonovits \cite[Theorem 2.2]{simonovits}
  implies (as a special case) that for a fixed number partition
  $\lambda \in \mathfrak{P}_m$, and for $n$ large enough, the largest graph
  on vertex set 
  $[n]$ that does not contain the clique graph corresponding to
  $\lambda$ is the complement of a clique
  graph. Phrased in our language, this means precisely that the
  antichain $\fib((\lambda, 1^{n-m}))$ in $\Pi_n$ has the Tur\'an
  property. By Corollary \ref{co:hook}, this antichain is not a blocker
  (unless $\lambda$ is a hook shape), so, by (1), its blocker ideal does
  not equal its vanishing ideal. Thus, we have another example showing
  that the converse of Proposition \ref{pr:turanideal} is false.}
\end{example}

\begin{example}{\rm 
It is easy to find antichains in $\Pi_n$ that do not satisfy the
Tur\'an property. One example is the antichain of any pair of atoms in
$\Pi_3$. 
 However, the only class of $\Sn$-invariant counterexamples that we
know of is the following. 

 Consider the partition
$\lambda=(2^r)$ for some $r$. Using Theorem \ref{th:blocker}, one
readily verifies that $\fib(\lambda)^* = \fib((r+1,1^{r-1}))$ in
$\Pi_{2r}$. Thus, the 
assertion that $\fib(\lambda)$ has the Tur\'an property  is
equivalent to the assertion that the smallest number of edges in any
graph on vertex set $[2r]$ that intersects every complete matching is
attained in an $(r+1)$-clique. However, a star (the graph containing
every possible edge from a single vertex) also intersects every
complete matching, and the star has fewer edges than the
$(r+1)$-clique if $r \geq 3$.}   
\end{example}

We end by describing a class of symmetric antichains with the Tur\'an property
which is not produced 
by Theorem \ref{pr:turanideal}.
Let $\Fq$ be the finite
field on $q$ elements, and consider the geometric lattice $L_q^n$ of all
subspaces of $\Fq^n$ ordered by inclusion. The analogue of orbit
arrangements would in this case be antichains that are invariant under
the action of $\GL (n,q)$, i.e.\ antichains $\A_k$ that contain every
subspace of a given dimension $k$.

Clearly, $\A_k^* = \A_{n-k+1}$.
The following proposition is therefore a reformulation of Theorem 3.5
in \cite[p.\ 87]{hirschfeld}, 
which says that a set of points in $PG(n,q)$ that intersects every
$k$-dimensional subspace has cardinality at least $1+q+\dots+q^{n-k}$.

\begin{proposition}
The antichain $\A_{k}\subset L_q^n$, which consists of all $k$-dimen\-sional
subspaces of  $\Fq^n$, has the Tur\'an property. 
\end{proposition}

\end{document}